\documentclass[12pt,a4paper]{amsart}
\usepackage{amsfonts,amsmath,amssymb}
\usepackage{hyperref}
\usepackage[all]{xy}

\newtheorem*{theorem*}{Theorem}
\newtheorem{lemma}{Lemma}[subsection]
\newtheorem{proposition}[lemma]{Proposition}
\newtheorem{remark}[lemma]{Remark}

\newtheorem{theorem}[lemma]{Theorem}
\newtheorem{definition}[lemma]{Definition}

\newtheorem{corollary}[lemma]{Corollary}

\baselineskip 18pt \textwidth 16cm  \theoremstyle{plain}

\newcommand{\cc}{\mathbb{C}}

\newcommand{\rr}{\mathbb{R}}

\newcommand{\Fre}{{Fr\'{e}chet \,}}

\newcommand{\cD}{{\mathcal{D}}}
\newcommand{\g}{{\mathfrak{g}}}

\newcommand{\gd}{\g^{\sigma}}

\begin{document}

\author{Eitan Sayag}

\address{Einstein Institute of Mathematics,
Edmond J. Safra Campus, Givat Ram, The Hebrew University of
Jerusalem, Jerusalem, 91904, Israel} \email{sayag@math.huji.ac.il}

\title{$(GL_{2n}(\cc),Sp_{2n}(\cc))$ is a Gelfand pair}
\date{\today}

\begin{abstract}
We prove that $(GL_{2n}(\cc),Sp_{2n}(\cc))$ is a Gelfand pair.
More precisely, we show that for an irreducible smooth admissible
Frechet representation $(\pi,E)$ of $GL_{2n}(\cc)$ the space of
continuous functionals $Hom_{Sp_{2n}(\cc)}(E,\cc)$ is at most one
dimensional. For this we show that any distribution on
$GL_{2n}(\cc)$ invariant with respect to the double action
$Sp_{2n}(\cc) \times Sp_{2n}(\cc)$ is transposition invariant.
Such a result was previously proven for $p$-adic fields by M.
Heumos and S. Rallis.

\end{abstract}

 \maketitle

\section{Introduction}\label{intro}

Let $F$ be a $p$-adic field (i.e. a local field of characteristic
zero). It was proven in \cite{MR1078382} that
$(GL_{2n}(F),Sp_{2n}(F))$ is a Gelfand pair. A simple proof, for the
case where $F$ is a finite field, appeared recently in
\cite{MR2290925}. In \cite{OS3} we used Proposition 3.1 of
\cite{MR2290925} to simplify the proof of the $p$-adic case (the
idea of proof, namely the usage of the Gelfand Kazhdan method
remains the same).

Our goal in this note is to transfer the result to the archimedean
case. This is the first step in extending the joint work of Omer
Offen and the author (see \cite{OS1}, \cite{OS2}, \cite{OS3}) to
the archimedean case.

 More precisely the main result of the present note is the following:
\begin{theorem*} [A]\label{main}
Let  $(\pi,E)$ be a irreducible admissible smooth Frechet
representation $(\pi,E)$ of $GL_{2n}(\cc)$. Then the space of
functionals $Hom_{Sp_{2n}(\cc)}(E,\cc)$ is at most one
dimensional.
\end{theorem*}

Following the standard technique of Gelfand and Kazhdan we first
show

\begin{theorem*} [B]\label{maindist}
Let $T$ be a distribution on $GL_{2n}(\cc)$ invariant with respect
to the double action $Sp_{2n}(\cc) \times Sp_{2n}(\cc)$. Then $T$
is transposition invariant.
\end{theorem*}

For the proof we will use an analogue of Gelfand-Kazhdan criterion
as proved in \cite{AGS} as well as the strategy of
\cite{AG2}-\cite{AG3}. Namely, we verify that the symmetric pair
$GL_{2n}(\cc),Sp_{2n}(\cc))$ is {\bf good} and that all its {\bf
descendants} are {\bf regular} (See section \ref{Prelim} for the
exact definitions of these terms).
We hope to cover the case $(GL_{2n}(\rr),Sp_{2n}(\rr))$ in the near
future.


We have made an attempt to make this note as self contained as
possible. Since our proof relies heavily on the notions and ideas
of \cite{AGS} and \cite{AG2} we summarized those in section
\ref{Prelim}. We introduce the notion of Gelfand pair and review
the Gelfand-Kazhdan distributional technique.

We then review some notions and results of \cite{AG2}, especially
the notions of symmetric pair, descendants of a symmetric pair,
good symmetric pair and regular symmetric pair.

The formal argument is given in the section \ref{formal} while the
key computation of descendants of our symmetric pair is done in
section \ref{descendance}.

\subsection*{Acknowledgements}
It is a pleasure to thank \textbf{Avraham Aizenbud} and
\textbf{Dmitry Gourevitch} for various discussions and for useful
advice.

\section{Notations and Preliminaries on Gelfand
Pairs}\label{Prelim}

In this section we review some notions and result concerning
Gelfand pairs. We follow the notions introduced in \cite{AGS} and
\cite{AG2}.

\subsection{Notation}\label{notation}

We let $F$ be a $p$-adic field, i.e. a local field of characteristic
zero. We denote algebraic $F$-varieties by bold face letters. Thus
${\bf X}$ is an $F$-variety and $X={\bf X}(F)$ is the set of its
$F$-points. Note that when $F$ is non-archimedean
 $X$ is an $\ell$-space (as in \cite{BZ}) while for archimedean $F$
 the set $X$ has a structure of smooth manifold.
We consider the space $S^{*}(X)$ of Schwartz distributions on $X$.
For $F$ non-archimedean these are just linear functionals on the
space $S(X)$ of locally constant functions with compact support on
the $\ell$-space $X$. When $F$ is archimedean we consider $X$ as a
Nash manifold and use the Schwartz space of \cite{AG1}. We note that
the main result of the paper, Theorem B, holds true when one
replaces $S^{*}(X)$ with $\cD(X)$ the space of all distributions on
$X$ (see Theorem 4.0.8 of \cite{AG2}).

The letters ${\bf G}, {\bf H}$ denotes reductive algebraic groups
defined over $F$. The groups of points of which will be denoted by
$G,H$ respectively. 

When we consider symmetric pairs we refer to triples ${\bf
X}=({\bf G},{\bf H},\theta)$ of reductive $F$-groups where
$\theta: {\bf G} \to {\bf G}$ is an involution and ${\bf H}={\bf
G}^{\theta}$ is the fixed point of $\theta$. The anti-involution
$\sigma(g)=\theta(g^{-1})$ and the symmetrization map
$s(g)=\sigma(g)g \in G^{\sigma}$ will play a key role.

For a symmetric pair ${\bf X}$ we will denote by $X(F)$ the pair
of groups $(G,H)$ of corresponding $F$-points of the groups ${\bf
G}, {\bf H}$.

If $A,B$ are matrices we will use the notation $diag(A,B)$ for the
block matrix whose blocks are $A$ and $B$.

For $n \geq 1$ we let $G_{n}=GL(n)$ and let $J_{2n}$ be the non
degenerate skew symmetric matrix of size $2n$ given by

\begin{equation*}
\begin{bmatrix}
0 & I_{n} \\
-I_{n} & 0
\end{bmatrix}
\end{equation*}

Clearly,
$$\theta_{2n}(g)=J_{2n}^{-1} {^tg} ^{-1} J_{2n}$$
is an involution on $G_{2n}$ and we have

$$H_{2n}=G_{2n}^{\theta_{2n}}=\{g \in G_{2n}:
^tgJ_{2n}g=J_{2n}\}=Sp(2n)$$

 We record the anti-involution $$\sigma_{2n}(g)=\theta_{n}(g^{-1})=J_{2n}^{-1} {^tg}
J_{2n}$$ and the symmetrization map $$s(g)=g J_{2n}^{-1} {^tg}
J_{2n}$$

\subsection{Distributional Criterion for Gelfand pairs} \label{GelPairs}
$ $\\

\begin{definition}
Let ${\bf G}$ be a reductive group. By an \textbf{admissible
representation of} $G$ we mean an admissible representation of
$G(F)$ if $F$ is non-archimedean (see \cite{BZ}) and admissible
smooth \Fre representation of $G(F)$ if $F$ is archimedean.
\end{definition}

We recall the following three notions of Gelfand pair.

\begin{definition}\label{GPs}
Let $H \subset G$ be a pair of reductive groups.
\begin{itemize}
\item We say that $(G,H)$ satisfy {\bf GP1} if for any irreducible
admissible representation $(\pi,E)$ of $G$ we have
$$\dim Hom_{H(F)}(E,\cc) \leq 1$$

\item We say that $(G,H)$ satisfy {\bf GP2} if for any irreducible
admissible representation $(\pi,E)$ of $G$ we have
$$\dim Hom_{H(F)}(E,\cc) \cdot \dim Hom_{H}(\widetilde{E},\cc)\leq
1$$

\item We say that $(G,H)$ satisfy {\bf GP3} if for any irreducible
{\bf unitary} representation $(\pi,\mathcal{H})$ of $G(F)$ on a
Hilbert space $\mathcal{H}$ we have
$$\dim Hom_{H(F)}(\mathcal{H}^{\infty},\cc) \leq 1.$$
\end{itemize}

\end{definition}
Property GP1 was established by Gelfand and Kazhdan in certain
$p$-adic cases (see \cite{GK}). Property GP2 was introduced in
\cite{Gross} in the $p$-adic setting. Property GP3 was studied
extensively by various authors under the name {\bf generalized
Gelfand pair} both in the real and $p$-adic settings (see e.g.
\cite{vD,Bos-vD}).

Clearly, $$GP1 \Longrightarrow GP2 \Longrightarrow GP3.$$

In \cite{AGS} we obtained the following version of a classical
theorem of Gelfand and Kazhdan.

\begin{theorem}\label{DistCrit}
Let $H \subset G$ be reductive groups and let $\tau$ be an
involutive anti-automorphism of $G$ and assume that $\tau(H)=H$.
Suppose $\tau(\xi)=\xi$ for all bi $H(F)$-invariant Schwartz
distributions $\xi$ on $G(F)$. Then $(G,H)$ satisfies GP2.
\end{theorem}

In the case $(G_{2n},H_{2n})$ at hand GP2 is equivalent to GP1 by
the following simple proposition (see \cite{AGS}).

\begin{proposition} \label{GP2GP1}
 Let $H \subset GL(n)$ be any transpose invariant subgroup.
 Then $GP1$ is equivalent to $GP2$ for the pair
$(\mathrm{GL}(n),H)$.
\end{proposition}

\subsection{Symmetric Pairs, GK pairs and Gelfand Pairs}

We now restrict to the case where $H \subset G$ is a symmetric
subgroup and more precisely $(G,H,\theta)$ is a symmetric pair.

\begin{definition}
A \textbf{symmetric pair} is a triple $(G,H,\theta)$ where ${\bf H}
\subset {\bf G}$ are reductive groups, and $\theta: {\bf G} \to {\bf
G}$ is an involution of ${\bf G}$ such that ${\bf H} = {\bf
G}^{\theta}$.
\end{definition}

In \cite{AG2} the notion of {\bf GK} symmetric pair is introduced.
It means that the condition of the above criterion (Theorem
\ref{DistCrit}) are met with $\tau$ being $\sigma$.

\begin{definition}
$(G,H,\theta)$ is called {\bf GK} pair if any $H(F) \times H(F)$
invariant distribution is $\sigma$ invariant.
\end{definition}

 In particular it follows from Theorem \ref{DistCrit} that
$$ GK \Longrightarrow GP2.$$

We now wish to describe the main tool to show that a given
symmetric pair is GK pair. In principle there is one obvious
obstruction for a pair to be GK pair. If there exist any closed
orbit $\Delta=HgH$ which is NOT $\sigma$ invariant then we can
find a distribution supported on that orbit which is not $\sigma$
invariant and this contradicts the condition. Thus an obvious
necessary condition is the following:

(*) for each closed orbit $HgH$ we have $\sigma(g) \in HgH$

Following \cite{AG2}, we call symmetric pairs satisfying condition
(*) {\bf good} symmetric pairs.

\begin{definition}
The symmetric pair $(G,H,\theta)$ is called good if for all closed
orbits $\Delta=HgH$ we have $\sigma(\Delta)=\Delta$.
\end{definition}

Thus, $$GK \Longrightarrow good$$

In \cite{AG2} it is conjectured that the converse is also true and
a partial result in that direction is proved.

To formulate the theorem we need two notions, that of {\bf
descendant} symmetric pair and that of {\bf regular} symmetric
pair.

The exact technical notion of {\bf regular} symmetric pairs will be
explained later and for the moment we note that the expectation is
that any symmetric pair is regular. Later we will discus this notion
and explain how we will verify that certain symmetric pairs are
regular.

The notion of {\bf descendant} pairs is central for the present
work. Given a symmetric pair $({\bf G},{\bf H},\theta)$ we
consider the centralizers in ${\bf G}$ and in ${\bf H}$ of certain
semi simple elements $x \in G(F)$. The elements $x$ that are
considered are in the image of the symmetrization map and in
particular satisfy $\sigma(x)=x$.

Being the centralizers of semi-simple elements these groups ${\bf
H_{x}} \subset {\bf G_{x}}$ are reductive, they are defined over $F$
and
$H_{x}$ is a
symmetric subgroup of $G_{x}$, namely the fixed point of
$\theta_{x}=\theta|_{G_{x}}$.

Let us now describe the elements $x \in G(F)$ that we consider.
Given $g \in G$ such that $\Delta(g)=HgH$ is closed it is known
that $x=s(g)$ is semi-simple and the stabilizer of $x$ in $H$ is
exactly the stabilizer of $g$ in $H \times H$ (for a proof see
Proposition 7.2.1 of \cite{AG2}). One can define a baby symmetric
pair $(G_{x},H_{x},\theta|_{G_{x}})$ by taking the centralizers of
$x$ in $G$ and $H$.

This is called a {\bf descendant} of $(G,H,\theta)$. Notice that the
set of descendants of a given symmetric pair depends on the field
$F$ over which the groups are considered.

\begin{remark}
It follows from standard finiteness results that the number of
non-conjugate descendants pairs of a given symmetric pair is always
finite (see corollary 1, p.107 of \cite{St}).
\end{remark}

We can now state one of the main results of \cite{AG2}:

\begin{theorem}\label{mainAG}
Suppose that the pair $(G,H,\theta)$ is {\bf good} and that all
its {\bf descendants} are {\bf regular}. Then $(G,H,\theta)$ is a
{\bf GK} pair.
\end{theorem}

\subsection{The notion of regularity}

To define the notion of regular symmetric pair we need several more
notations. Given an algebraic representation $V$ of $G$ we say that
$v \in V(F)$ is nilpotent if zero is in the closure of its orbit
$G(F)v$. We denote by $\Gamma(V)$ the set of nilpotent elements.

Consider the $Q(V)=(V/V^{G})(F)$ and the subspace $R(V) \subset
Q(V)$ of regular element (i.e. those that are not nilpotent).
Notice that since $G$ is reductive we can realize $Q(V)$ as a
subspace of $V(F)$. Thus,
$$Q(V)=V/V^{G}$$ and $$R(V)=Q(V)-\Gamma(V)$$ where $\Gamma(V)$ is the set of
nilpotent elements for the action of $G$ on $V$.

Let $(G,H,\theta)$ be a symmetric space.

We now formulate two properties $Dist_{gen}(g)$ and
$Dist_{all}(g)$.

$$Dist_{gen}(g,V) \Longleftrightarrow S^{*}(R(V))^{H(F)} \subset
S^{*}(R(V))^{Ad(g)}$$ and

$$Dist_{all}(g,V) \Longleftrightarrow S^{*}(Q(V))^{H(F)} \subset
S^{*}(Q(V))^{Ad(g)}$$

We remind that $S^{*}(Q(V))$ is the usual space of Schwartz
distributions on a vector space, while for the definition of
$S^{*}(R(V))$ we use the theory developed in \cite{AG1}.

Now consider $\g:=LieG$ and the adjoint action of $H(F)$ on the
vector space $\gd:=\{a \in \g | \theta(a)=-a\}$.

We expect that for some elements $g \in G(F)$ the fact that each
$H(F)$-invariant  distribution on the regular part $R(V) \subset
Q(V)$ is $Ad(g)$ invariant will imply that the same is true for
distributions on the entire vector space $Q(V)$.

We now specify the elements $g \in G(F)$ that will play a role in
the definition of regularity. We deviate here slightly from
\cite{AG2}. Let $Z=Z(G)$ be the center of $G$.

An element $g \in G$ is called {\bf relevant} if $s(g):=g \sigma(g)
\in Z(G)$ and $H$ has index at most 2 in the group $H<g>$ generated
in $Aut(G)$ by $Ad(H)$ and $Ad(g)$.

We denote by $$rel(G,H,\theta)=\{g \in G: s(g) \in Z(G),
[H<g>:Ad(H)] \leq 2\}$$

the set of relevant elements.

 We can now give the definition of regularity.

\begin{definition}
The pair $(G,H)$ is called {\bf regular} if for the adjoint action
of $H(F)$ on the vector space $\gd$
$$Dist_{gen}(g,\gd)
\Longrightarrow Dist_{all}(g,\gd)$$ for all relevant elements $g
\in rel(G,H,\theta)$.

\end{definition}

\begin{remark}
In \cite{AG2} the authors use the notion of admissible elements.
These satisfy a further condition which is implied by
$Dist_{gen}(g,\gd)$. Namely, if $Dist_{gen}(g,\gd)$ holds then
each closed $H(F)$ orbit on $\gd$ is preserved by $Ad(g)$. This is
exactly the condition {\bf $Ad(g)|_{\gd}$ is $H$-admissible}. Thus
the definition given here for regularity is identical to that of
\cite{AG2}.
\end{remark}

We have the following trivial consequence of the definitions.
\begin{proposition}\label{regularityCrit}
We have the following:

(i) the product of regular pairs is a regular pair.

(ii)  Let $(G,H,\theta)$ be a symmetric pair. $Z=Z(G)$ the center of
$G$. Suppose that
$$\{g: \sigma(g)g \in Z\} \subset ZH.$$ Then $(G,H,\theta)$ is
regular.

\end{proposition}

\begin{proof}
The first point is proposition 7.4.4 of \cite{AG2} (We emphasize
that this is a direct consequence of the definitions and the
Proposition 3.1.5 of \cite{AGS}).

For the second point we note that $rel(G,H,\theta)$ is by
assumptions a subset of $ZH$ and hence the condition
$Dist_{all}(g)$ is obviously satisfied.
\end{proof}

\section{${\bf X_{n}}(F):=(G_{2n},H_{2n},\theta)$ satisfies GK when $F \ne \rr$}\label{formal}

In this section we prove Theorem A and Theorem B.

For the proof of Theorem A it is enough to show that the pair
$(G_{2n},H_{2n},\theta)$ is GK.

Indeed, by Theorem \ref{DistCrit} it satisfies GP2 and hence, by
Proposition \ref{GP2GP1} it satisfies GP1. This is the content of
Theorem A.

Observe that Theorem B is equivalent to the assertion that the
pair $(G_{2n},H_{2n},\theta)$ is GK because the anti-involution
$\sigma_{2n}$ is conjugate to transposition via the element
$J_{2n}$ of $Sp_{2n}$.

 To show that $X_{n}$ is GK we will verify that it is {\bf good} and compute its
descendants.

When $F=\cc$ we show that all these descendants are regular. Thus by
Theorem \ref{mainAG} the pair $X_{n}(\cc)$ is GK.

\begin{remark}
It can be shown that $X_{n}(F)$ is regular also for $F$
non-archimedean. We don't need this fact here since our proof that
$X_{n}(F)$ is GK will be based on the original Gelfand-Kazhdan
criterion that requires the verification that {\bf all} the orbits
$\Delta=HgH$ are $\sigma$-invariant. See below.
\end{remark}

\subsection{The pair $X_{n}(F)$ is good}

We first verify that all orbits are transpose invariants.

\begin{proposition}\label{transpose invariant}
Let $F$ be an arbitrary field. Let $g \in GL_{2n}(F)$. Then $g^{t}
\in HgH$ where $H=Sp_{2n}(F)$.
\end{proposition}
The proof depends on the following lemma (which is proposition 3.1
of \cite{MR2290925}):

\begin{lemma}\label{for good}
Let $F$ be an arbitrary field. For  $x \in GL_{n}(F)$ define
$d(x)=diag(x,I_{n})$. The map $Ad(G_{n})x \to H_{2n}d(x)H_{2n}$ is a
bijection between the set of conjugacy classes in $GL_{n}(F)$ and
the set of orbits of $Sp_{2n}(F) \times Sp_{2n}(F)$ in $GL_{2n}(F)$.
In particular, given $g \in GL_{2n}(F)$ the double coset
$H_{2n}gH_{2n}$ contains an element of the form $d(x)=diag(x,I_{n})$
where $x \in GL_{n}(F)$.
\end{lemma}

 {\bf Proof of proposition
\ref{transpose invariant}}

Since $H$ is transpose invariant, the statement is equivalent to
 $$HgH=(HgH)^{t}$$ Thus we can choose whatever representative we
 like in $HgH$. We take $g=diag(x,I_{n})$
 as in the lemma. But for
 some $k \in GL(n)$ we have $kxk^{-1}=x^{t}$ since in $GL(n)$ any
 matrix is conjugate to its transposed. Taking $h=diag(k,{^tk}^{-1})$
 which belongs to $H_{2n}$ we have $$h \bullet
 diag(x,I_{n}) \bullet h^{-1}=diag(^tx,I_{n})$$ as needed.

\begin{corollary}
Let $F$ be a $p$-adic field. Then the pair $(GL_{2n}(F),Sp_{2n}(F))$
is good.
\end{corollary}
Indeed, by the above proposition and $J_{2n} \in H_{2n}$ we see that
for each orbit $\Delta=HgH$ we have $$\sigma(g)=J^tgJ \in HgH$$

\begin{remark}
For the case $F=\cc$, it is shown in \cite{AG2} (see corollary
7.1.7) that any connected symmetric pair (i.e. $G/H$ is connected)
is {\bf good}.
\end{remark}

\begin{remark}
The fact that the pair $X_{n}(F)$ is good for any local field $F$
can be deduced from the {\bf cohomological criterion} given in
corollary 7.1.5 of \cite{AG2}. Indeed, all descendants are products
of symplectic groups (as shown in section \ref{descendance} of this
note) and these are known to be cohomologically trivial (\cite{Sp}).
In a sense the approach of \cite{MR1078382} is along these lines.
\end{remark}

Since all orbits of $H_{2n}(F) \times H_{2n}(F)$ on $G_{2n}(F)$ are
$\sigma$-invariant one can use, in the case that $F$ is a $p$-adic
field the work of Gelfand-Kazhdan to deduce

\begin{theorem}
Let $F$ be a non-archimedean field. Then the pair
$(G_{2n},H_{2n})$ is a GK and hence a Gelfand pair.
\end{theorem}

The complete argument can be extracted from page 6 of \cite{OS3}
but it is really an immediate consequence of the theory of
Gelfand-Kazhdan as presented by in \cite{BZ}.

\begin{remark}
Note that when $F$ is non-archimedean the condition of regularity
is not required in our case since {\bf all the orbits} for the
pair $X_{n}$ are $\sigma$-invariant.
\end{remark}

\subsection{The descendants of the pair $(GL_{2n}(F),Sp_{2n}(F))$}

We have

\begin{proposition}\label{alldesc}

All descendants of $X_{n}(F)$ are products of symmetric pairs of
the form $X_{m}(E)$ for some $m \leq n$ and some $E/F$ finite
extension.
\end{proposition}
For the proof see the fourth section where the proposition is
formulated in slightly different form.

\subsection{Regularity of $X_{n}$ and its descendants, $F=\cc$}

We now restrict to the case $F=\cc$.

\begin{corollary}
Let $F=\cc$. The symmetric pair $X_{n}(F)$ and all its descendants
are regular.
\end{corollary}

\begin{proof}
By proposition \ref{alldesc}, each descendant is a product of
$X_{m}(E)$ for some field extension $E/F$ and $m<n$. Thus, by
clause (i) of proposition \ref{regularityCrit}, it is enough to
prove that $X_{m}(E)$ is regular for any $E/F$ and any $m$.

Let $g$ be a relevant element. Then $g \sigma(g) \in Z(G_{2m})$ and
since $F=\cc$ we can find a scalar matrix $\alpha \in Z(G_{2m})$ and
$g_{0} \in H_{2m}$ such that $g=\alpha g_{0}$. Thus we search for
$\alpha \in Z(G_{2m})$ such that $\alpha^{-1}g \in H$. We claim that
any $\alpha$ solving the equation

\begin{equation}\label{eq}
\alpha \sigma(\alpha)=z
\end{equation}
Will work. Indeed,
$$\alpha^{-1} g \sigma(g) \sigma(\alpha)^{-1}=1$$ (we use
$z=g\sigma(g)$) and thus $s(\alpha^{-1}g)=1$ showing that
$\alpha^{-1}g \in H$. We thus need to solve the equation (\ref{eq})
or $$\alpha J^{-1} \alpha J=z$$ This can be done by hands after
writing $z=diag(t,...,t)$ and taking $\alpha=diag(r,r,..,r)$ with
$r=t^{1/2}$.
By clause (ii) of proposition \ref{regularityCrit} we are done.
\end{proof}

\section{The descendants of the pair $(GL(V),Sp(V))$}\label{descendance}

Let $F$ be arbitrary local field. Let $(V,\omega)$ be a symplectic
space over $F$. Let $J:V \to V^{*}$ be the isomorphism induced by
$\omega$. $J$ is given by $J(v)(u)=\omega(v,u)$. Clearly, $J^{*}=-J$
where $J^{*}$ is the adjoint operator to $J$.

Consider the involution $\theta(g)=J^{-1} g^{\tau} J$ where
$g^{\tau}$ is the inverse of the adjoint of $g:V \to V$. Note that
the isomorphism $\tau:GL(V) \to GL(V^{*})$ given by
$\tau(g)=(g^{*})^{-1}$ correspond in matrix notations to transpose
inverse.

The fixed point group of $\theta$ is the symplectic group
$$Sp(V,J)=\{g \in GL(V): Jg=g^{\tau}J\}$$
Clearly, this is the group of automorphisms of the symplectic space
$(V,\omega)$.
The symmetric pair
$(GL(V),Sp(V,J),\theta)$ will be denoted $X_{V,J}$.

The aim of this section is to verify the following
\begin{theorem}\label{descendant}
Let $F$ be arbitrary local field. Let $X_{V,J}$ be the symmetric
pair as above. Then all descendants of $X_{V,J}$ are products of
pairs of the form $X_{W,I}$ where $(W,I)$ is a symplectic space
over a finite extension $E/F$ with $dim_{E}(W) \leq dim_{F}(V)$.
\end{theorem}

For the proof of the main result of the note we need the above
only for the case $F=\cc$. In that case it can be proven by brute
force. One can treat arbitrary local field using the method of
\cite{SpSt} as we now show.

We recall the formula for the centralizer of a semi-simple
transformation. Remind that for a semi-simple element $x \in GL(V)$
the minimal polynomial $P=min(x,F)$ is a product of distinct
$F$-irreducible polynomials. Moreover, if $P=\prod_{i=1}^{n} P_{i}$
then

\begin{equation}\label{jordan}
V= \oplus V_{i}
\end{equation}
 with
$V_{i}=Ker(P_{i}(x))$. Note that the adjoint operator $x^{t}$ yields
a decomposition $V^{*}=\oplus V^{*}_{i}$ where
$V^{*}_{i}=Ker(P_{i}(x^{*}))$. Here $x^{*}:V^{*} \to V^{*}$ is
defined by $$x^{*}(\ell)(v)=\ell(xv).$$ We have
$x^{\tau}=(x^{*})^{-1}$.

Notice that each $V_{i}$ is a vector space over the field
$E_{i}=F[t]/(P_{i}(t))$. Also, an element that commutes with $x$
must preserve the direct sum decomposition $V= \oplus V_{i}$ and
clearly commutes with the action of each $E_{i}$ on $V_{i}$.
Thus,

$$GL(V)_{x}=\prod_{i=1}^{n} GL(V_{i},E_{i})$$

The content of proposition 2.8 (page E-89) in \cite{SpSt} is an
extension of this well known fact to describe the centralizers of
elements of classical groups. We need a related result, namely a
description of centralizers of elements that satisfy $x=\sigma(x)$
(as opposed to $x=\theta(x)$). Let us present the result in the case
where $V$ is endowed with a symplectic form $\omega$ or equivalently
with $J:V \to V^{*}$ which is an isomorphism and anti-self dual. We
remind that $\theta(x)=J^{-1}(x^{*})^{-1}J$ and
$\sigma(x)=J^{-1}x^{*}J$.

We will require two simple lemmas.

\begin{lemma}
Suppose $x=\sigma(x)$. Then the decomposition (\ref{jordan})
$V=\oplus V_{i}$ is orthogonal with respect to $\omega$.
\end{lemma}

\begin{proof}
First note that the condition $Jx=x^{*}J$ is equivalent to
$\omega(xv,u)=\omega(v,xu)$. Let $v \in Ker(P(x))$ and $u \in
ker(Q(x))$ with $P,Q$ distinct factors and hence prime to each
other. Find $a(x),b(x)$ with
$$v=a(x)P(x)v+b(x)Q(x)v.$$ We have
$$\omega(v,u)=\omega(b(x)Q(x)v,u)=\omega(b(x)v,Q(x)u)=0$$
\end{proof}

\begin{lemma}
Let $x \in GL(V)$ a semi-simple element with $x=\sigma(x)$. Let
$P(x)=\prod_{i=1}^{n}P_{i}(x)$ be the decomposition of the minimal
polynomial of $x$ to irreducible factors. Then $J:V \to V^{*}$ maps
each $V_{i}=ker(P_{i}(x))$ to $V_{i}^{*}=Ker(P_{i}(x^{*}))$.
Moreover, $J$ is $E_{i}$ linear on each $V_{i}$ and $J|_{V_{i}}$
defines a non-degenerate form on the space $V_{i}$.

\end{lemma}

\begin{proof}
If $v \in V_{i}$ then $P_{i}(x)v=0$. Now since $Jx=x^{*}J$ we get
$P_{i}(x^{*})Jv=JP_{i}(x)v=0$ and hence $Jv \in V^{*}_{i}$. The
$E_{i}$ linearity is checked using the same identity. We thus see
that $J= \oplus J_{i}$ with $J_{i}=J|_{V_{i}}$ and since $J$ is
invertible, each of the maps $J_{i}$ is invertible. Since $J^{*}=-J$
and $J^{*}=\oplus J^{*}_{i}$ we see that each of the forms $J_{i}$
is skew symmetric and invertible thus defines a symplectic space
$(V_{i},J_{i})$ over $E_{i}$.
\end{proof}

\begin{proof}[Proof of Theorem \ref{descendant}]

Let $g \in GL(V)$ be an element which is semi-simple with respect to
the $Sp(V) \times Sp(V)$ action. The corresponding descendant is
defined via the element $x=s(g)=\sigma(g)g$ which is semi-simple and
belongs to $GL(V)^{\sigma}$. So let $x \in GL(V)^{\sigma}$ be a
semi-simple element. We claim that $(G_{x},H_{x})$ is a product of
symmetric pairs as above. By the comment above

$$G_{x}
=\prod_{i=1}^{n} GL(V_{i},E_{i})$$

Let now $h \in H=Sp(V,J)$ be an element in $H_{x} \subset G_{x}$.
Since $h \in H$ we have $h^{*}J=Jh$. Since $h \in G_{x}$ we can
write everything in blocks corresponding to the decompositions of
$V$ and $V^{*}$.  Using the previous lemma we obtain that
$h=diag(h_{1},...,h_{r})$ with each $h_{i}:V_{i} \to V_{i}$ satisfy
$h^{*}_{i}=J_{i}h_{i}J_{i}^{-1}$. Thus,

$$H_{x}
=\prod_{i=1}^{n} Sp((V_{i},J_{i})/E_{i})$$

This is exactly the claim.

\end{proof}

\end{document}